\documentclass[12pt]{article}
\usepackage{amssymb, amscd}
\newcommand{\Ha}{\mathbb{H}}
\newcommand{\C}{\mathbb{C}}

\newcommand{\Q}{\mathbb{Q}}

\newcommand{\Z}{{\mathbb{Z}}}

\newcommand{{\1}}{{\bf 1}}

\newcommand{\co}{{\mathcal{O}}}

\newcommand{\ce}{{\mathcal{E}}}

\newcommand{\ra}{\rightarrow}

\newcommand{\lra}{\longrightarrow}

\newcommand{\vp}{\varphi}

 \newcommand{\IM}{\mbox{Im }}
\newcommand{\Aut}{\mbox{Aut}}

\newcommand{\Ker}{\mbox{\rm ker}}
\newcommand{\End}{\mbox{End}}
\newcommand{\Hom}{\mbox{Hom}}
\newcommand{\Pic}{\mbox{Pic}}
\newcommand{\Div}{\mbox{Div}}

\def \D{ \displaystyle}

\begin{document}
\def\flexleftarrow#1{\mathtop{\makebox[#1]{\leftarrowfill}}}
\begin{center}
{\Large \bf Abelian varieties with group action}\\[2ex]
by\\[1ex]
H. Lange \footnote{Supported by CONACYT 27962-E} and S. Recillas \footnote{Supported by CONACYT 
  27962-E
and DAAD}
\end{center}
%\title[Abelian varieties with group action]{Abelian varieties with group action}
%
%\author{Herbert Lange}
%\author{Sevin Recillas}
%\address{H. Lange\\Mathematisches Institut\\
%              Universit\"at Erlangen-N\"urnberg\\
%              Bismarckstra\ss e $1\frac{ 1}{2}$\\
%             D-$91054$ Erlangen\\
%              Germany}
%\email{lange@mi.uni-erlangen.de}
%\address{S. Recillas\\Instituto de Matem\'aticas\\
%                      UNAM Campus Morelia\\
%                      Morelia, Mich., 58190\\
%                     M\'exico\\
%                      and\\
%                      CIMAT\\ 
%                      Callejon de Jalisco, Gto., 58000\\
%                      M\'exico\email{sevin@matmor.unam.mx} 
%\thanks{Supported by CONACYT 27962-E and DAAD}
%\keywords{abelian variety, Jacobian variety, Prym variety}
%\subjclass[2000]{Primary: 14H40; Secondary: 14k02}

\vspace{0.5cm}
\begin{abstract}
\noindent
A finite group $G$ acting on an abelian variety $A$ induces a decomposition of $A$ up
to isogeny. In this paper we investigate several aspects of this decomposition. We apply 
the results to the decomposition of the Jacobian variety of a smooth projective curve with $G$-action. 
The aim is to express the decomposition in terms of the $G$-action of the curve. As examples of the results we work
out the decompositions of Jacobians with an action of the symmetric group of degree 4, 
the alternating group of degree 5, the dihedral groups 
of order $2p$ and $4p$, and the quaternion group.
\end{abstract}

\vspace{1cm}
{\large \bf Introduction} 
\vspace{0.5cm}

\noindent
Let $G$ be a finite group acting on a smooth projective curve $X$. This induces an
action of $G$ on the Jacobian $JX$ of $X$ and thus a decomposition  of $JX$ up to isogeny. The
most prominent example of such a situation is the case of the group $G \simeq \Z/2 \Z$ of two elements.
Let $\pi : X \ra Y = X/G$ denote the canonical quotient map. The first to notice that $JX$ is isogenous
to the product $JY \times P(X/Y)$ was probably Wirtinger (see [W]). The abelian variety $P(X/Y)$ was
later called by Mumford the Prym variety of the covering $\pi$, because of its relation to certain Prym
differentials.
\smallskip

\noindent
The group action of some other special groups were studied by Recillas who in [Re] 
decomposed the Jacobian of the Galois covering of a simple trigonal cover, Ries who in [Ri] 
decomposed the Jacobian of a $D_p$-gonal curve,
 Recillas-Rodr\'{\i}guez
($S_3$ in [RR1] and $S_4$ in [RR2]), Donagi-Markman ($S_3, S_4$ and $W D_4$ in [DM]),
Carocca-Recillas-Rodr\'{\i}guez (the dihedral groups in [CRR]), and S\'{a}nchez-Arg\'{a}ez ($A_5$ in [SA2]). There
are also some general results: First of all there is the classical paper of Chevalley-Weil (see [CW])
which computes how many times an irreducible complex representation appears in the tangent space $T_0
JX$ in terms of the fixed points under elements of $G$  and its stabilizer subgroups. Finally Donagi in
[D] and M\'{e}rindol in [M] found the decomposition of $JX$ in the case of a finite group $G$ all of
whose irreducible $\Q$-representations are absolutely irreducible.
\medskip

\noindent
It is the aim of the present paper to prove some general results on abelian varieties with action by a
finite group $G$ and deduce from them  a method to decompose the Jacobian $JX$ of a curve $X$ in terms of the
G-action of the curve. It will be applied in several examples.\\

\noindent
Our point of view for doing this is the following:
Let $G$ be a finite group acting on an abelian variety $A$. The action induces an algebra homomorphism
$\rho : \Q [G] \ra \End_{\Q} (A)$. We define $\IM (\alpha) := \IM (\rho (m \alpha)) \subseteq A$ for any
$\alpha \in \Q [G]$ where $m$ is some positive integer such that $m \alpha \in \Z [G]$. Let
$$
1 = e_1 + \ldots + e_r
$$
denote the decomposition of $1 \in \Q [G]$ into the sum of central orthogonal idempotents $e_i$ of $\Q [G]$
corresponding to the simple $\Q$-algebras associated to $\Q [G]$. 
If $A_i = \IM (e_i)$, then it is obvious that the $A_i$ are $G$-stable abelian subvarieties of $A$ with
$\Hom_G (A_i, A_j) =0$ for $i \not= j$ such that the addition map induces an isogeny
$$
A_1 \times \ldots \times A_r \ra A,
$$
called the {\it isotypical decomposition } of $A$. The essential point here is that one can explicitly
describe the elements $e_i$ in terms of the complex characters of $G$, which allows to actually compute
the isotypical decomposition in many cases. In Section 2 we decompose the $G$-abelian subvarieties
$A_i$ further. In fact, if
$$
e_i = f_{i_1} + \ldots + f_{i_{n_i}}
$$
is any decomposition of $e_i$ into a sum of orthogonal primitive idempotents, this
 gives rise to an isogeny decomposition
$$
B_{i_1} \times \ldots \times B_{i_{n_i}} \ra A_i
$$
with $B_{i_j} = \IM (f_{i_j})$. One can show that
$B_{i_1} , \ldots , B_{i_{n_i}}$ are isogenous to each other, so that finally we obtain an isogeny
decomposition
$$
A \sim B^{n_1}_{1} \times \ldots \times B^{n_r}_{r} \eqno(1)
$$
Here the abelian varieties $B_i$ correspond to the different irreducible $\Q$-representations of
$G$.
This generalizes a result of M\'{e}rindol [M] who deduced the same decomposition in the special case that
all irreducible $\Q$-representations of $G$ are absolutely irreducible. 
\medskip

\noindent
In Section 3 we apply these results to the case of a Jacobian variety with $G$-action. Let $X$ be a
smooth projective curve of genus $\geq 2$ over the field of complex numbers and let $G$ be a finite
group acting on $X$ with quotient $\pi : X \ra Y  := X / G$. Let $W_1 , \ldots , W_r$ denote the
irreducible $\Q$-representations of $G$, where we assume that $W_1$ is the trivial representation. Then
if 
$n_i$ denotes the dimension of $W_i$ as a $D_i$-left vector space where $D_i = \End_G (W_i)$, then the
isogeny decomposition (1) of the Jacobian $JX$ reads 
$$
JX \sim    \pi^{\ast} J Y \times B^{n_2}_2 \times \ldots \times B_r^{n_r}. \eqno(2)
$$
The first problem is of course to see which of the abelian varieties $B_i$ are nonzero. This certainly depends
on the curve $X$ and not only on the group $G$. However, we show in Theorem 3.1 that if the genus of
$Y$ is $\geq 2$, then $B_i \not= 0$ for all $i$.
\\
The second problem is to describe the abelian variety $B_i$ in terms of $G$-action of the the curve $X$ and its
quotients, that is for example in terms of the Prym varieties of the intermediate coverings of $\pi : X \ra Y$.
Recall that if $M \subset N \subset G$ are subgroups with intermediate covering $\vp : X_M \ra X_N$ where $X_M = X /
M$ and $X_N = X / N$, then the {\it Prym variety} $P(X_M/X_N)$ of $\vp$ is defined as the complement of
$\vp^{\ast} JX_N$ in $JX_M$ with respect to the canonical polarization and we have
$$
JX_M \sim    \vp^{\ast} JX_N \times P(X_M / X_N).
$$
It is a priori not clear that the isogeny decomposition of $P(X_M / X_N)$ underlies the isogeny
decomposition (2), since the abelian varieties $B_i$ are not simple in general. But Proposition 3.7
 says that this is the case. In fact, we show that there are non negative integers $s_i$ such that
$$
P(X_M / X_N) \sim    B^{s_2}_2 \times \ldots \times B^{s_r}_r.
$$
This can be used to compute the isogeny decomposition of $JX$ in many cases. In Section 4 we work out
the examples $G = S_4, A_5$, the dihedral groups $D_p$ and $D_{2p}$ for an odd prime $p$ and the quaternion 
group $Q_8$.\\

Although some of the results are valid for abelian varieties over an arbitrary algebraically closed fields,
we need for the essential theorems the description of an abelian variety as a complex torus. So for
the sake of simplicity we assume the base field to be the field $\C$ of complex numbers.\\

\medskip
\noindent
{\bf Acknowlegements:} The first author would like to thank Gabriele Nebe for teaching him some aspects of 
the theory of $\Q$-representations. The second author wants to thank A. S\'{a}nchez-Arg\'{a}ez for helpful
conversations. This work was carried out during a visit of the first author  to Morelia and of 
the second author 
to Erlangen under the DAAD-UNAM exchange program. We also thank the ICTP, Trieste, where the final version of
the paper was written.
\vspace{1cm}

{\large \bf 1. The isotypical decomposition}

\vspace{0.5cm}
\noindent
Let A be an abelian variety over the field $\C$ of complex numbers and 
let $G$ denote a finite group acting on $A$. The action induces an
algebra homomorphism
$$
\rho : \Q [G] \ra \End_{\Q} (A)
$$
where $\Q [G]$ denotes the group algebra of $G$ over $\Q$. 
For the sake of simplicity we often consider an element $\alpha \in \Q [G]$ also as an element of
$\End_{\Q} (A)$ i.e. we write $\alpha$ instead of $\rho (\alpha)$.
 
If $\alpha$ is any
element in $\Q [G]$, we define
$$
\IM (\alpha) := \IM (\rho ( m \alpha)) \subseteq A
$$
where $m$ is some positive integer such that $m \alpha \in
\Z [G]$. $\IM (\alpha)$ is
an abelian subvariety of $A$, which certainly does not depend on the chosen
integer $m$. In order to obtain proper abelian subvarieties we have to choose
suitable elements $\alpha$ of $\Q[G]$.
\medskip

\noindent
For this recall that $\Q [G]$ is a semisimple $\Q$-algebra of finite
dimension. As such it admits a unique decomposition
$$
\Q [G] = Q_1 \times \ldots \times Q_r \eqno(1.1)
$$
with simple $\Q$-algebras $Q_1, \ldots , Q_r$. Let
$$
1 = e_1 + \ldots + e_r \eqno(1.2)
$$
denote the decomposition of the unit element 1 given by (1.1). The elements
$e_i$ are the unit elements of the algebras $Q_i$ and form, when considered as
elements of $\Q [G]$, a set of orthogonal idempotents contained in the
center of $\Q [G]$.
\medskip

\noindent
{\bf Proposition 1.1:} {\it Let} $A_i = \IM (e_i)$ {\it for $i = 1, \ldots , r$}.
\begin{itemize}
\item[(a)] {\it $A_i$  is a $G$-stable abelian subvariety of $A$ with} $\Hom_G(A_i, A_j) =0$  {\it for $i \not= j$.
\item[(b)] The addition map induces an isogeny}
$$
\mu : A_1 \times \ldots \times A_r \ra A.
$$
\end{itemize}
This decomposition is called the {\it isotypical decomposition} of $A$. It is
unique up to permutation, since the idempotents $e_i$ are uniquely determined by
$\Q[G]$.
\medskip

\noindent
{\it Proof:} (a) Let $n_i$ be a positive integer  such that $n_i e_i \in \Z [G]$. The idempotent
$e_i$ being in the center of $\Q [G]$, we obtain
$$
g A_i = g (n_i e_i )A = n_i e_i g A \subseteq n_i e_i A = A_i.
$$
The second assertion follows from the fact that the $Q_i$ are simple $\Q$-algebras with $Q_i \not\simeq Q_j$ 
for $i \not= j$.
Assertion (b) is an immediate consequence of (1.2). \hfill $\square$
\medskip

\noindent
In order to obtain the isotypical decomposition of $A$, one proceeds as follows:
The action of $G$ induces a commutative diagram
$$
\begin{array}{ccc}
&& \End (T_0 A)\\
& \stackrel{\tilde{\rho}}{\nearrow} &\\
\Q [G] && \uparrow \rho_a \\
& \stackrel{\searrow}{\scriptstyle \rho} &\\
&& \End_{\Q} (A)
\end{array}
\eqno(1.3)
$$
where $T_0 A$ denotes the tangent space of $A$ in $0$ and $\rho_a$ is the
analytic representation of $\End_{\Q} (A)$. So $G$ acts compatibly on the
tangent space $T_0A$.
\smallskip

\noindent
Denote by $\Lambda$ the lattice in $T_0A$ such that $A \simeq T_0A/ \Lambda$. Defining $U_i = \IM (\tilde{\rho} (e_i))$, we obtain
$$
A_i = U_i / \Lambda \cap U_i
$$
for $i=1, \ldots, r$ and the isogeny $\mu : A_1 \times \ldots \times A_r \ra A$ induces a $G$-decomposition
$$
T_0 A = U_1 \oplus \ldots \oplus U_{r}.
$$
If $W_i$ denotes the irreducible $\Q$-representation of $G$ corresponding to the $\Q$-algebra $Q_i$ (in fact, $W_i$ is a minimal left ideal of
$Q_i$), we call $U_i$ {\it the $G$-subspace of $T_0A$ associated to } $W_i$. The idea for computing these subspaces is to express the
idempotents $e_i$ in terms of the representations of $G$.\\
This can be done as follows:

Let $V_1, \ldots , V_s$ denote the irreducible
$\C$-representations  of the group $G$ with corresponding characters $\chi_1,
\ldots , \chi_s$. Then
$$
p_j := \frac{\deg \chi_j}{|G|} \,\, \sum_{g \in G} \overline{\chi_j (g)} g \in \C [G]
$$
is a projector of $\C [G]$ with image the sum of all sub-representations of $\C
[G]$ isomorphic to $V_j$.\\
Since the values $\chi_j (g)$ are roots of unity, the field
$$
K_j = \Q (\chi_j (g),  g \in G)
$$
is cyclotomic of degree say $d_j$ over $\Q$. In particular $K_j \, | \, \Q$ is a
Galois extension. The conjugate representations of $V_j$, i.e. the
representations associated to the conjugate characters $\sigma \chi_j \,\,\, (
\sigma \in$ Gal $(K_j | \Q))$ are again irreducible over $\C$ and non 
isomorphic to
each other. Hence they are among the $V_i$'s. Let $V_{j_1} = V_j, V_{j_2} ,
\ldots, V_{j_{d_j}}$ be the different irreducible representations conjugate to
$V_j$. Their projectors $p_{j_{\nu}}$ are pairwise orthogonal and their sum
$p_{j_1} + \ldots + p_{j_{d_j}}$ is again a projector. It is a projector of $\Q
[G]$, since
$$
\begin{array}{rcl}
p_{j_1} + \ldots + p_{j_{d_j}} &=& \D{ \frac{ \deg \chi_j}{|G|} \,\, \sum_{g \in
G}
(\overline{\chi_{j_1} (g)} + \ldots + \overline{\chi_{j_{d_j}} (g)}) \cdot g}\\[3ex]
&=& \D{ \frac{ \deg \chi_j}{|G|} \,\, \sum_{g \in G} tr_{K_j|\Q} (\chi_j (g)) \cdot g \,\,\,\,
\in \Q [G]}.
\end{array}
$$
Here we use the fact that all characters $\chi_{j_{\nu}}$ are of the same
degree. The image of the projector $p_{j_1} + \ldots + p_{j_{d_j}}$ in
$\C [G]$ consists of the sum of all subrepresentations of $\C [G]$ isomorphic
to one of the $V_{j_{\nu}}$. 
But this is just $Q_j$, the image of the idempotent $e_j$. Hence we have
$$
e_{j} = p_{j_1} + \ldots + p_{j_{d_j}} = \frac{\deg \chi_j}{|G|} \,\, \sum_{g
\in G} tr_{K_j |\Q} (\chi_j (g)) \cdot g. \eqno(1.4)
$$

\medskip
\noindent

\vspace{1cm}
{\large \bf 2. Decomposition of the isotypical components}
\vspace{0.5cm}

\noindent
Let $G$ be a finite group acting on an abelian variety $A$ with isotypical decomposition $\mu : A_1
\times \ldots \times A_r \ra A$. The components $A_i$ correspond one to one to the irreducible $\Q-$representations 
$W_i$ of $G$. In order to avoid indices, let $A^{0}$ denote one of them. Let $W$ denote the corresponding 
irreducible $\Q$-representation and $Q$ the associated simple $\Q-$subalgebra of $\Q [G]$ with unit element $e_W$ 
given by equations (1.1) and (1.4). According to Schur's Lemma
$$ 
D := \End_G(W)
$$ 
is a skew-field of finite dimension over $\Q$. Then $W$ is a finite dimensional left-vector space over $D$ and we have 
$$
 Q \simeq \End_D (W). \eqno(2.1) 
$$
Let $ n:= \dim_D (W) $. Choosing a $D-$basis of $W$ identifies $\End_D (W)$ with the 
matrix ring $M(n \times n, D)$. 
Let $ E_i \in M(n \times n, D)$ denote the matrix with entry $1$ at the place $(i,i)$ and $0$ elsewhere
and let $p_i \in Q$ denote the element corresponding to it under the isomorphism (2.1). Then $ \{ p_1, \ldots , p_n \} $
is a set of primitive orthogonal idempotents in $ Q \subset \Q [G]$ with
$$
e_W = p_1 + \cdots + p_n.  \eqno(2.2)
$$
Consider the abelian subvarieties $ B_i = \IM (p_i)$ of $A^0 $ as defined at the beginning of Section 1. 
Equation (2.2) immediately implies that the addition map
$$
\mu : B_1 \times \ldots \times B_n \ra A^0 
$$
is an isogeny. The idempotents $p_1, \ldots , p_n$ are certainly conjugate to each other in the 
algebra $Q$, since 
$$
E_i = P_{ij} E_j P_{ij}^{-1}
$$
where $P_{ij} \in M(n \times n, D) $ denotes the permutation matrix corresponding to the transposition
$(i,j)$. This implies that $B_1, \ldots , B_n$ are pairwise isogenous. Hence we obtain 
\noindent
\medskip

\noindent
{\bf Proposition 2.1:} {\it Let $A^{0} \not= 0$ be the isotypical component of an
abelian variety $A$ with $G$-action associated to an irreducible $\Q-$representation $W$ and} $n:= \dim_D(W)$ 
{\it with} $D := \End_G(W)$. {\it Then there is an abelian subvariety $B$ of $A^{0}$ such that $A^0$ is isogenous to $B^n$:
$$
A^{0}  \sim B^{n}.
$$.}

\noindent
Combining Propositions 1.1 and 2.1 we obtain
\medskip

\noindent
{\bf Theorem 2.2:} {\it Let $G$ be a finite group acting on an abelian variety $A$. Let $W_1, \ldots ,
W_r$ denote the irreducible $\Q$-representations of $G$ and} $n_i:= \dim_{D_i}(W_i)$ 
{\it with} $D_i := \End_G(W_i)$ {\it for $i = 1,
\ldots , r$. Then there are abelian subvarieties $B_1, \ldots , B_r$ of $A$ and an isogeny}
$$
A \sim B^{n_1}_1 \times \ldots \times B^{n_r}_r.
$$
We call this decomposition an {\it isogeny decomposition of $A$ with respect to} $G$. Note that 
the abelian subvarieties $B_i$ of $A$ are only determined up to isogeny. 
 Note moreover that the abelian varieties $B_i$ might be zero. For example if one starts with an
isotypical component $A_i$, its isogeny decomposition is just $A_i \sim B_i^{n_i}$. So even 
then one gets
a proper decomposition if only $ n_i > 1$.\\

\noindent
As elements of $\Q [G]$ the idempotents $p_i$ of (2.2) are of the form
$$
p_i = \sum_{g \in G} r_g^i \ g
$$
with rational numbers $r_g^i$. In special cases one can be more precise. Suppose $W$ is an absolutely 
irreducible representation of $G$, that is the extension $ W \otimes_{\Q} \C$ of $W$ is an irreducible 
$\C-$representation. Hence $D = \Q$ and $W$ is of dimension $n$ over $\Q$. In this case $W$ 
admits up to a positive constant a uniquely determined $G-$invariant scalar product (see [Se]). Fix one
of these and denote it by $(\ \ ,\ \, \, )$, let $\{ w_1, \ldots , w_n \}$ be a basis of $W$, 
orthogonal with respect to $(\ \ ,\, \,)$, and define
$$
p_{w_i} := \frac{n}{|G| \cdot \| w_i \|^2} \sum_{g \in G} (w_i,gw_i) \ g.
$$
It was communicated to us by Gabriele Nebe that Schur's character relations can be translated into terms of idempotents as 
follows\\

\noindent
{\bf Proposition 2.3:} {\it $p_{w_1}, \ldots , p_{w_n}$ are orthogonal idempotents in $\Q [G]$ satisfying
$$
p_{w_1}+ \ldots + p_{w_n} = e_W.
$$}

{\it Proof:} For any $ g \in G$ let $(r_{ij}(g))_{i,j}$ denote the matrix of $g$ with respect to the basis 
$\{ w_1, \ldots , w_n \}$ of $W$. Schur's character relations are (see [Se], Chapter 2, Corollary 3 of Proposition 4)
$$
\sum_{g \in G} r_{ij}(g^{-1})r_{kl}(g) = \frac{|G|}{n}\  \delta_{il} \ \delta_{jk}
$$
where $\delta_{ij}$ is equal to $1$ if $i=j$ and $0$ otherwise. Now
$$ 
(w_i,gw_i) = (g^{-1}w_i,w_i) = (\sum_{j=1}^{n} r_{ij}(g^{-1})w_j,w_i) = r_{ii}(g^{-1}) \|w_i\|^2
$$
and thus applying the character relations we obtain the following matrix equation
$$
\left( \left(\sum_{g \in G} (w_i,gw_i)r_{kl}(g)\right)_{k,l}\right) = \|w_i\|^2 
\left(\left(\sum_{g \in G} r_{ii}(g^{-1})r_{kl}(g)\right)_{k,l}\right)
= \|w_i\|^2 \frac{|G|}{n} \cdot E_i
$$
where $E_i$ denotes the matrix of above. Since the matrices $E_i$ are idempotents in the algebra $M(n \times n, \Q)$ 
with $E_1 + \cdots + E_n = I_n$, 
this implies the assertion.           \hfill $\square$\\

\noindent
{\bf Remark 2.4:} 
Proposition 2.3 is valid in slightly more generality, namely for all irreducible $\Q-$representations which remain 
irreducible over the field of {\it real} numbers. In fact, in [M] M\'erindol gives another 
proof of Proposition 2.3, which is valid
for all $\Q-$representations which admit only a $1-$dimensional space of $G-$invariant symmetric bilinear form. Note that a 
$\Q-$representation is of this kind if and only if it remains irreducible when 
extended to a real representation.
 
\vspace{1cm}

{\large \bf 3. Jacobians with group action}
\vspace{0.5cm}  

\noindent
In this section we apply the previous results to the Jacobian varieties  of
algebraic curves. 
The curve $X$ is assumed to be smooth, irreducible, projective, defined
over $\C$ and of genus $g_X \geq 2$. We will denote by $JX = H^0 (X, \omega_X)^{\vee} / H_1 (X,
\Z)$  its Jacobian variety, which is canonically isomorphic to
$\Pic^0 X = H^1 (X, \co_X)/ H^1 (X, \Z)$.

In the first two sections we considered group actions on unpolarized abelian varieties. Here it will 
be important that the Jacobians are polarized abelian varieties. In fact, the intersection product
on $H_1 (X, \Z)$ induces a canonical principal polarization $H = c_1(L)$ on $JX$. Here $L$ denotes the 
line bundle on $JX$ defined by the theta divisor. Its first Chern class $H$ can be considered as 
a positive definite hermitian form on the tangent space $T_0JX$
whose imaginary part $E$ is integer valued on the lattice $\Lambda =  H_1 (X, \Z)$. 
The polarization $H$ induces an isomorphism
$\phi_{H}: JX \stackrel{\simeq}{\lra} \widehat{JX}$ of JX onto its dual abelian variety $\widehat{JX}$. 
An automorphism $g : JX
\ra JX$ in the sequel always means a polarization preserving one meaning that $g^*L \approx L$
(algebraically equivalent) or equivalently that $\phi_{g^*H} = \phi_H$. 
Note that $g \in \Aut (X)$ preserves the polarization $H$ if and only if $H(gv, gw) =
H(v,w)$ for all $v,w \in T_0JX$.\\

\noindent
If a curve $X$ has a non trivial group of automorphisms, then by
Torelli's theorem we have:
$$
\Aut (X) = \left\{
\begin{array}{ll}
\Aut (J X) & \mbox{if $X$ is hyperelliptic},\\
\Aut (JX)/\{\pm 1\} & \mbox{if not}.
\end{array}
\right.
$$
One knows that $|\Aut (X) \, | \, \leq  84 (g_X -1)$ and also that {\it any} finite
group $G$ is a group of automorphisms of a compact Riemann surface (see
[H]).\\

\noindent
Let us now consider a curve $X$ with group action $G \subset \Aut X$,
so $G$ acts on the Jacobian $JX$. 
According to Section 1 the tangent space of $JX$ at $0$ decomposes as follows 
$$
 T_0 JX = H^1(X, \co_X)
 =  U_1 \oplus \ldots \oplus U_r,
$$ 
where $U_i$ denote the $G$-subspaces of $T_0 JX$ associated to the irreducible $\Q$-representations 
$W_i$ of $G$. If $\Lambda = H^1(X,\Z)$, then 
$$
A_i = U_i / \Lambda \cap U_i     \eqno(3.1)
$$ 
are abelian subvarieties of $JX$ and from Proposition 1.1 we know that
the addition  map
$$
\mu : A_1 \times \ldots \times A_r \ra JX     \eqno(3.2)
$$
is an isogeny.
Recall from Proposition 2.1 that if  $n_i = \dim_{D_i} W_i$ where $D_i = \End_G W_i$, then $A_i$ is
isogenous to the $n_i$-fold product of an abelian subvariety $B_i$
with itself, $A_i \sim  B_i^{n_i}$. So the isogeny  decomposition of $JX$ is:
$$
J X \sim B_1^{n_1} \times \ldots \times B^{n_r}_r . \eqno(3.3)
$$
Our aim is to deduce an isogeny decomposition of $JX$ from the lattice of subgroups $H$ of $G$ and the induced actions.
In the sequel we always denote by $V_1$ the trivial and by $V_2, \ldots , V_s$ the 
nontrivial irreducible $\C$-representations of $G$ and similarly by $W_1$ the trivial and by
$W_2, \ldots , W_r$ the nontrivial $\Q$-representations of $G$. Then $\C \otimes_{\Q} W_i$ is the direct sum
$$
\C \otimes_{\Q} W_i = V_{j_1} \oplus \ldots \oplus V_{j_{s_i}}
$$ 
of some of the $V_j$'s (to
each $j$ there is a unique $i$), with $s_i = d_i m_i$ where $d_i$ is as at the end of Section 1 and 
$m_i$ is the Schur index of the representation $W_i$.\\

\noindent
First we
consider the question which of the abelian varieties $B_i$ in (3.1) are nonzero. The next theorem implies
that if the genus of the quotient
$Y = X/G$ is $ \geq 2$, then $B_i \not= 0$ for all $i$.\\

\noindent
{\bf Theorem 3.1:} {\it Given a Galois covering of smooth projective curves $\pi : X \ra Y$,
with group $G$, all irreducible $\Q$-representations of $G$ appear in the
isotypical decomposition of $JX$ if $g_Y \geq 2$.}\\

\noindent
For this we need the following lemma which reflects the fact that $\pi_{\ast} \co_X$ may be considered
as a sheafified version of the regular representation of $G$:\\

\noindent
{\bf Lemma 3.2}: {\it Let $X$ be a smooth projective curve with an action of $G$. Denote by $Y = X/G$ 
the quotient curve and by $\pi: X \lra Y$ the covering map. Then
$$
\pi_{\ast} \co_X = \co_{Y} \oplus (V_2 \otimes_{\C} \ce_{V_2}) \oplus
\ldots \oplus (V_s \otimes_{\C} \ce_{V_s})    \eqno(3.4)
$$
where the $\ce_{V_i}$ are locally free $\co_{Y}   $-modules such that $rk
\ce_{V_i} = \dim V_i$ and the representations $V_i$ are as above.}
\bigskip

{\it Proof}: $\pi_{\ast} \co_X$ is torsion free and thus locally free of rank $n = \#G$, since $X$ is a curve. 
The action of $G$ on the curve $X$ induces an action on the sheaf $\pi_{\ast} \co_X$  over the trivial 
action of $G$ on $Y$. Let $\C(X)$ respectively $\C(Y)$ denote the function fields of $X$ and $Y$.
For any $G$-invariant open affine set $ U = Spec(R) \subset Y$ trivializing
$\pi_{\ast} \co_X$, the action of $G = Gal(\C(X)/\C(Y))$ on $\pi_{\ast} \co_X(U) = R^n \subset \C(Y)^n$  
is the restriction of the action of $G$ on $\C(Y)^n = \C(X)$, which is the regular representation of $G$.
Since this is independent of $U$, we get if we sheafify
$$
\pi_{\ast} \co_X = \co_Y \oplus (V_2 \otimes_{\C} \ce_{V_2}) \oplus
\ldots \oplus (V_s \otimes_{\C} \ce_{V_s}) 
$$
where the $\ce_{V_i}$ are locally free $\co_Y$-modules such that $rk\ce_{V_i} = dim V_i$. 
\hfill $\square$\\

\noindent
{\it Proof of Theorem 3.1:} From (3.4) we deduce that
$$
H^1 (X, \co_X) \cong  H^1 (Y, \pi_{\ast} \co_X) \cong  H^1 (Y, \co_Y)
\oplus \bigoplus^s_{j=2} H^1 (Y, \ce_{V_j})^{\dim V_j} .
$$
Hence it suffices to show that $h^1 (Y, \ce_{V_j})>0$ for all $i \geq 2$, since then
all representations $W_i$ will appear in the isotypical decomposition of 
$ T_0 J X$.\\
\noindent
By Grothendieck duality (see [K]) one knows that $(\pi_{\ast}
\co_X)^{\vee} \cong  \pi_{\ast} \omega_{X/Y}$, where $\omega_{X/Y}$ 
denotes the relative dualizing sheaf. Since $\pi$ is a finite
morphism, we know that $\pi_{\ast} \omega_{X/Y}$ is nef (see [V]). Hence all
quotients of $(\pi_{\ast} \co_X)^{\vee}$ have non negative degree.
Since 
$$
(\pi_{\ast} \co_X)^{\vee} \cong  \co_Y \oplus \D{ \bigoplus^s_{j=2}} (V_j
\otimes_{\C} \ce_{V_j})^{\vee}
$$ 
it follows that $\deg (\ce_{V_j}^{\vee})
\geq 0$. \\
On the other hand from (3.4) and $H^0(Y, \co_Y) \cong \C \cong H^0 (X, \co_X) \cong  H^0 (Y, \pi_{\ast} \co_X)$ we deduce
$h^0 (Y, \ce_{V_j}) = \{0\}$ for $j = 2, \ldots, s$.  
So from Riemann-Roch  we obtain
$$
h^1 (Y, \ce_{V_j}) = - \deg \ce_{V_j} + \dim V_j (g_Y -1) > 0
$$
if $g_Y \geq 2$. \hfill $\square$\\

\noindent
We may assume that each $B_i$ is an   abelian subvariety  of $JX$. In fact, $B_i$ may be considered as
 the image of an
endomorphism of $JX$. We give now   correspondences on $X$ which induce such
endomorphisms.

\noindent
Given our Galois cover of curves $\pi : X \ra Y$ with Galois group
$G$,
for each $g \in G$ we consider its graph:
$$
\Gamma_g = \{ (x, gx); \, x \in X \} \subset X \times X.
$$
Given any element  
 $\D{ p = \sum^{}_{g \in G} n_g g \in \Z [ G ]}$ 
 we consider the divisor
$$
D := \sum_{g \in G} n_g 
\Gamma_g \in \mbox{Div} (X \times X)
$$
which defines a correspondence on $X$, namely:
$$
X  \ra  \D{ \mbox{Div}^{\sum n_g} }(X), \quad
x \mapsto  \D{ \sum_{g \in G} n_g g(x)}.
$$
\noindent
Let
$$
\gamma_D : JX \ra JX.
$$
be the induced endomorphism. In other words  $\gamma_D = \rho (\sum_{g \in G} n_g g)$.\\

\noindent
An example of this is when $n_g =1$ for all $g \in G$. In this case we
have $\D{ \sum\nolimits^{}_{g \in G} g x = \pi^{\ast} (\pi x)}$, so it follows
that $\IM \gamma_D  = \pi ^{\ast} JY$ for $ D = \D{ \sum\nolimits^{}_{g \in G} 
\Gamma_g}$.\\

\noindent
For every irreducible $\Q-$representation $W_i$ of $G$ let $e_{W_i} = p_{i,1} + \cdots + p_{i,{n_i}} $
denote the decomposition (2.2) of $e_{W_i}$ into a sum of primitive orthogonal idempotents of the simple 
subalgebra $Q_i \subset \Q [G]$ associated to $W_i$. Fix positive integers $m_i$ such that  
$$
 q_i := m_i p_{i,1} = \sum_{g \in G} n_g^i g  \in \Z [G]. \eqno (3.5)
$$
Then if
$$
B_i := \IM (p_{i,1}),
$$ 
we conclude\\

\noindent
{\bf Proposition 3.3:}  \, \,{\it $B_i = \IM \gamma_{\Sigma n^i_g \Gamma_g }$.}\\

\noindent
Let us now consider a morphism of curves $f : X \ra Y$ of degree $d$. Let $H = H_X$ denote the canonical 
principal polarization of $JX$. We
define the {\it complement} of $f^{\ast} JY$ on $JX$ with respect to $H$ as
$$
P := P (X/Y) := \Ker( \widehat{f^{\ast}} \circ \phi_H )^{\circ} = \Ker (N m
f)^{\circ}. 
$$
where $Nmf: JX \ra JY$ denotes the norm map. 
We call $P$ the {\it Prym variety}
of $f$ (even if it is not a Prym variety in the usual sense, i.e. if the principal polarization does not restrict
to a multiple of a principal polarization on $P$).
       Denote by $H_P$ the induced polarization  on $P$ and 
observe that the induced polarization on $f^{\ast} JY$ is $d \cdot {H_Y}$
(see [LB]). The      map 
$$
\sigma : (JY,d \cdot {H_Y}) \times (P(X/Y), H_P) \ra (JX, H),
\quad (a,b) \ra
f^{\ast}  a + b
$$ is an isogeny of polarized abelian varieties. 
Applying a simple induction this yields:\\

\noindent
{\bf Proposition 3.4:} {\it For a tower of curves
$X \stackrel{f_1}{\ra} X_1 \stackrel{f_2}{\ra} X_2 \ra \ldots \stackrel{f_r}{\ra}
X_r$, denote $g_i = f_i \circ \ldots \circ f_1$. Then the decomposition}
$$
T_0 JX = T_0 P(X/X_1) \oplus (d g^{\ast}_1)_{0} T_0 P (X_1 /X_2) \oplus \ldots \oplus
(d g^{\ast}_{r-1})_0 T_0 P (X_{r-1} /X_r) \oplus (d g^{\ast}_r)_0 T_0 J X_r 
$$
{\it is $H$-orthogonal}. \hfill $\square$ \\

\noindent
We can say more about the role of the trivial representation of $G$ in the isogeny decomposition (3.3).
Since the action of the group $G$ on $JX$ is induced by the action of $G$ on the curve $X$,  $G$ preserves the 
canonical polarization $H$. Hence the isotypical decomposition
$$
T_0 JX \cong U_1 \oplus \ldots \oplus U_r
$$
is $H$-orthogonal.
We always assume that $U_1$ is the $G$-subspace of $T_0 JX$ associated to the trivial representation.
On the other hand if $\pi : X \ra Y = X/G$ as above, by Proposition 3.4 the decomposition
$$
T_0 JX \cong (d \pi^{\ast})_0 (T_0 JY) \oplus T_0 P(X/Y)
$$
is $H$-orthogonal. Since we have seen above that
$$
(d \pi^{\ast} )_0 (T_0 JY) \cong U_1 
$$
and since both decompositions of $T_0JX$ are $H$-orthogonal, this implies
$$
T_0 P(X/Y) \cong \bigoplus^r_{i=2} U_i  .
$$
We call this argument {\it orthogonal cancellation.}\\

\noindent
With the notation of (3.1) and (3.2) it follows from this that
$$
\pi^{\ast} JY = A_1 \quad \mbox{and} \quad P(X/Y) \sim    A_2 \times \ldots
\times A_r
$$
with $A_1 = \IM \D{ \sum_{g \in G} g}$ and $n_1 = 1$.\\

\noindent
One can also get information about the abelian subvarieties $B_i \subset JX$ for
the non trivial representations of $G$.
For  this let $M$ be a subgroup of $G$, $X_M = X/M$ and $\vp : X \ra X_M$ the quotient
map.
\noindent
For $g \in G$ we consider the divisor
$$
\overline{\Gamma}_g := (\vp \times \vp)_{\ast} \Gamma_g \in \,\, \mbox{Div}  (X_M \times X_M) .
$$
Since  $\vp(x) = \vp (y)$ for $x,y \in X$ if and only if there is an  $ h \in M$ such that $y = hx, \quad
\overline{\Gamma}_g$ induces the  map:
$$
\overline{g} :X_M  \longrightarrow  \Div^{|M|} (X_M)\qquad \qquad 
z  \mapsto          \overline{g}(z) =  \sum_{h \in M} \vp (ghx)
$$
where $x \in X$ is such that $\vp(x) =z$. This definition does not depend on the choice of $x$.\\
\noindent
More generally, given a divisor $D = \D{ \sum\nolimits_{g \in G} n_g \Gamma_g \in \Div (X \times X)}$,
consider the divisor $\overline{D} = \D{ \sum\nolimits_{g \in G} n_g \overline{\Gamma}_g
\in \Div (X_M \times X_M)}$.  
Denote by $\gamma_D \in \End (JX)$ and
$\gamma_{\overline{D}} \in \End (JX_M)$ the induced endomorphisms. Then we have\\

\noindent
{\bf Proposition 3.5:} {\it If for each $g \in G$ one has $n_g = n_{gh} =
n_{hg}$ for all $h \in M$, then the following diagram is commutative:
$$
\begin{array}{ccc}
JX & \stackrel{\gamma_D}{\longrightarrow} & JX\\
{\scriptstyle N m \vp} \downarrow \qquad && \qquad \downarrow {\scriptstyle N m \vp}\\
JX_M & \stackrel{\frac{1}{|M|} \gamma_{\bar{D}}}{\longrightarrow} & JX_M\\
{\scriptstyle \vp^{\ast}} \downarrow \quad && \quad \downarrow {\scriptstyle \vp^{\ast}}\\
JX & \stackrel{\gamma_D}{\longrightarrow} & JX
\end{array}
$$
Moreover }$\IM \gamma_D = \vp^{\ast} \IM \gamma_{\overline{D}} \subset
\vp^{\ast} JX_M$.
\medskip

\noindent
{\it Proof.} Given $x \in X$, let $z = \vp (x) \in X_M$. The commutativity of the upper
square follows from:
$$
\sum_{g \in G} n_g \overline{g}(z) = \sum_{g \in G} n_g \sum_{h \in M} \vp (ghx) = \sum_{h
\in M} \sum_{ g \in G} n_g \vp ( g h x) = |M| \cdot \sum_{g \in G} n_g \varphi (g x).
$$
Here the last equality is a consequence of the conditions satisfied by the coefficients. A
similar computation shows that the lower square is commutative. 
The last assertion follows from the equation
$$
\vp^{\ast} \bar{D} (z) = \sum_{k
\in M} \sum_{ g \in G} \sum_{ h \in M} n_g k g h (x) = |M|^2 \cdot D(x) \hspace{3.5cm} \square
$$ \\

\noindent
{\bf Corollary 3.6:} {\it Let $q_i := m_i p_{i,1} = \sum_{g \in G} n_g^i g  \in \Z [G]$ 
for $i \geq 2$ be as in {\rm (3.5)}
and denote $ D_i = \sum_{g \in G} n^i_g \Gamma_g$. If $q_i h = h q_i = q_i$ for all $h \in M$, then}
$$
B_i = \IM (p_{i,1}) = \IM \gamma_{D_i}  \subset \vp^{\ast} P(X_M/Y) .
$$
{\it Proof:} By Proposition 3.3 we have $\IM \gamma_{D_i} = \IM (q_i)$. Since $q_i 
\cdot \sum_{g \in G} g  =0$ in $\Q [G]$, it follows that $\IM \gamma_{D_i}
\subset A_2 \times \ldots \times A_r$, hence by orthogonal cancellation $\IM
\gamma_{D_i} \subset P(X/Y)$.
The hypotheses on $q_i$ just mean that the assumptions of Proposition 3.5 are fullfilled and we obtain 
$\IM \gamma_{D_i} \subset \vp^{\ast} J X_M$, hence $\IM \gamma_{D_i} \subset
(\vp^{\ast} J X_M \cap P(X/Y))^0 = \vp^{\ast} P(X_M/Y)$ . \hfill $\square$
\medskip

\noindent
Observe that we have $B_i \sim    \IM \gamma_{\bar{D}_i} \subset P(X_M/Y).
$
So a natural question arises here: 
Given the decomposition (3.3) 
$$ 
J X \sim    JY \times B^{n_2}_2 \times \ldots \times B^{n_r}_r.
$$
Is there an underlying decomposition for $P(X_M / Y)$?
Observe that this is obvious if the abelian varieties $B_i$ are simple, which
 need not be the case. However we can show the following theorem which gives an answer
to this question in the slightly more general case of an intermediate
cover.\\
\medskip

\noindent
{\bf Proposition 3.7:} {\it  Given a Galois cover $ \pi : X \ra Y = X/G$ and subgroups $M
\subset N \subset G$, let $\psi : X_M \ra X_N$ be the associated cover where $X_N = X/N$ and
$X_M = X/M$. Let $JX \sim     JY \times B^{n_2}_2 \times \ldots \times B^{n_r}_r$
be the decomposition of above. Then there exist non negative
integers $s_i , \,\, i=2, \ldots , r$ such that}
$$
P(X_M / X_N) \sim    B^{s_2}_2 \times \ldots \times B^{s_r}_r.
$$
\noindent
{\it Proof: } Given the subgroups $M \subset N \subset G$ we have the tower of curves:
$$
\begin{array}{ccc}
X &&\\
& \searrow \vp & \\
&&X_M \\
\pi \downarrow && \downarrow \psi\\
&&X_N \\
&\swarrow \eta &\\
Y &&
\end{array}
$$
First we claim that it suffices to show that 
$$
 JX_M \sim JY \times B_2^{t_2} \times \cdots \times B_r^{t_r}  \eqno(3.6)
$$
with nonnegative integers $t_2, \ldots , t_r$ for every subgroup $M$ of $G$.\\
\noindent
To see this note that then we have the $H-$orthogonal decompositions
$$
T_0 JX_M \cong T_0 JY \oplus T_0 (B_2^{t_2}) \oplus \cdots \oplus T_0 (B_r^{t_r})
$$
$$
T_0 JX_N \cong T_0 JY \oplus T_0 (B_2^{t'_2}) \oplus \cdots \oplus T_0 (B_r^{t'_r})
$$ 
with nonnegative integers $ t'_i \leq t_i$ for $ i=1, \ldots , r$. On the other hand
$$
T_0 JX_M \cong T_0 JX_N \oplus T_0 P(X_M/X_N).
$$
Hence by orthogonal cancellation we obtain
$$
T_0 (B_2^{t'_2}) \oplus \cdots \oplus 
T_0 (B_r^{t'_r}) \oplus T_0 P(X_M/X_N) \cong T_0 (B_2^{t_2}) \oplus \cdots \oplus T_0 (B_r^{t_r})
$$
This gives the assertion, with $s_i = t_i - t'_i$.\\
\noindent
For the proof of (3.6) note that from invariant theory we know that
$$
(d \vp^{\ast})_0  T_0 JX_M = (T_0 JX)^M. 
$$
Hence $T_0 JX_M$ corresponds to the induced representation $ind^G_{1_M}$ of the trivial representation of 
$M$ in $G$ considered as a subrepresentation of $\Q [G]$. Since $JY$ and the $B_i$ correspond 
to the trivial and nontrivial irreducible $\Q-$representations of $G$, it follows that $T_0 JX_M$ is a direct sum of
$T_0 JY$ and some $T_0 B_i$'s, where $T_0 JY$ occurs only once since it corresponds 
to the trivial representation of $G$. This implies the assertion. \hfill $\square$\\ 

\noindent
For the following proposition let the notation be as in Proposition 3.7 and recall that $W_2, \ldots , W_r$ 
denote the nontrivial irreducible $\Q$-representations of $G$. Every irreducible $\C$-representation 
is contained in exactly one irreducible $\Q$-representation. Hence for every $W_i$ we may choose 
an irreducible $\C$-representation    
associated to $W_i$ and denote it by $V_i$.   
\medskip

\noindent
{\bf Proposition 3.8:} {\it Suppose for some $ i \geq 2$}\\
$$
\dim V^M_j - \dim V^N_j = \left\{
\begin{array}{cll}
1 && j = i \\
&\mbox{if}&\\
0 && j \neq i.
\end{array}
\right.
$$\\
{\it Then} \hspace{4cm} $P(X_M / X_N) \sim    B_i$.\\

\noindent
{\it Proof:} By assumption $V_i$ contains a fixed point of $M$ which is not 
a fixed point of $N$. Hence also the conjugate representations of $V_i$
over $\Q$ contain such a fixed point. This implies that also $W_i$ contains a fixed point, say
$w_1$, of M which is not a fixed point of $N$. Choosing $w_1$ as the first vector of the 
$D_i$-basis ${w_1, \ldots, w_{n_i} }$ of $W_i$ (where $D_i = \End_G (W_i)$),
all $h \in M$ are of the form $ h = \pmatrix{ 1 & 0 \cr 0 & h' \cr }$ with respect to this basis. 
The idempotent $p = p_{i,1} = \pmatrix{ 1 & 0 \cr 0 & 0 \cr }$ satisfies 
$ hp = ph = p$ for all $h \in M$. According to Corollary 3.6 the abelian subvariety $B_i = \IM (p)$ 
is contained in $\vp^{\ast} P(X_N/Y)$. Since $B_i$ is contained in $\vp^{\ast}JX_M$
and not in $ (\psi \vp)^{\ast}JX_N$, this implies that 
$$
B_i \subset P(X_M / X_N).
$$
In a similar way the assumption $ V_j^M = V_j^N$ for $ j \neq i$ implies that none of the $B_j, \,\,\, 
j \neq i$ is contained in $P(X_M / X_N)$, that is $s_j = 0$ for $j \neq i$ 
in the decomposition of Proposition 3.7.\\
\noindent
It remains to show that $s_i = 1$. But $s_i \geq 2$ would imply that $W_i$ admits two fixed points 
of $M$ which are linearly independent over $D_i$. 
This contradicts the assumption $\dim V_i^M = \dim V_i^N + 1$. \hfill $\square$\\

\noindent
A method for computing the decomposition of the Jacobian $JX$ of a curve
$X$ with group action was given in [SA] in the case $G = A_5$. Proposition 3.8
 allows us to deduce a generalization of this method which will be
useful for computing other examples. The starting point is the following
lemma. 
Here for any subgroup $M \subset G$ we denote by $\chi_{ind^G_{1_M}}$ the
character of the induced representation $ind^G_{1_M}$ of the trivial
representation of $M$ in $G$.\\

\noindent
{\bf Lemma 3.9:} {\it Let $G$ be a finite group, $V_1 (= \C), V_2, \ldots , V_s$ its
irreducible complex representations. Let $M \subset N \subset G$ be subgroups. If we
have for some $1 \leq l \leq s$ :}
$$
\chi_{ind^G_{1_M}} - \chi_{ind^G_{1_N}} = \sum^l_{i=1} \dim_{\C} V^M_i \cdot
\chi_{V_i},               \eqno(3.7)
$$
{\it then} 
$$
V^N_j = \left\{
\begin{array}{cll}
0 && 1 \leq j \leq l\\
&\mbox{if}&\\
V^M_j && l+1 \leq j \leq s.
\end{array}
\right.
$$
\medskip

\noindent
In other words, if we write
$$
ind^G_{1_M} = ind^G_{1_N} \oplus V_2^{\oplus b_2} \oplus \cdots \oplus V_l^{\oplus b_l}
$$
with $V_1, \ldots, V_l$ irreducible complex representations and positive integers $b_i$,
then the assumption (3.7) implies that none of the $V_i$'s appears in the representation 
$ind^G_{1_N}$.\\

\noindent
{\it Proof.} By Frobenius reciprocity:
$$
\begin{array}{lcl}
\dim_{\C} V^N_j &=& <\chi_{V_j} , \chi_{ind^G_{1_N}}>\\[2ex]
&=& < \chi_{V_j} , \chi_{ind^G_{1_M}} - \D{ \sum^l_{i=1} \dim V^M_i \cdot
\chi_{V_i}}>\\[2ex]
&=& \dim V^M_j -  \D{ \sum^l_{i=1} \dim V^M_i < \chi_{V_j},
\chi_{V_i}>}\\[4ex]
&=& \left\{
\begin{array}{ccl}
0& \mbox{if}& 1 \leq j \leq l\\
\dim V^M_j & \mbox{if}& l < j \leq s. \hspace*{4cm}\square
\end{array}
\right.
\end{array}
$$\\

\noindent
{\bf Corollary 3.10:} {\it Let $X$ be a curve with $G$ as a group
of automorphisms and $M \subset N \subset G$ subgroups. Let $B_j, W_j$ and $V_j$ be as at the beginning of this section and let 
$V_{j,1} = V_j, V_{j,2}, \ldots , V_{j,d_j}$ denote the nonisomorphic irreducible $\C$-representations
in $W_j \otimes_{\Q} \C$.\\
Assume that for some $i \geq 2$ we have}\\
\medskip
(i)  $\dim V^M_i =1$ {\it and}\\ 
\noindent
(ii) $\chi_{ind^G_{1_M}} - \chi_{ind^G_{1_N}} = \sum^{d_i}_{k=1} \chi_{V_{i,k}}$.\\
\medskip
\noindent
{\it Then}
$$
P ( X_M /X_N) \sim B_i.
$$\\

\noindent
{\it Proof:} Observe first that $\dim (V_{i,k})^M = \dim V^M_i =1$, for all \, $k, \,\, 1 \leq k \leq d_i$, 
since the representation $V_{i,k}$ is conjugate to $V_i$ over $\Q$. Hence by 
Lemma 3.9 we have:
$$
\dim (V_{j,k})^N  = \left\{
\begin{array}{cll}
0 && j = i \\
&\mbox{if}&\\
\dim (V_{j,k})^M && j \neq i.
\end{array}
\right.
$$
So Proposition 3.8 implies $P(X_M/ X_N) \sim    B_i$.
\hfill $\square$\\

\vspace{1cm}
{\large \bf 4. Examples}
\vspace{0.5cm}
\newline
Let $X$ be a smooth projective curve with $G$-action, $\pi :X \ra Y = X/G$ the
corresponding Galois covering. If $W$ denotes an irreducible $\Q$-representation, we will denote by
$A_W$ the isotypical component of $JX$ corresponding to $W$. If $D = \End_G (W)$ and $ n = \dim_D (W)$ we know 
by Proposition 2.1 that there is an abelian subvariety $B_W$ of $A_W$ such that 
$$
A_W \sim B_W^n.
$$
Note that $A_W$ is uniquely determined whereas $B_W$ is not in general. The aim of this section is to apply the above 
results to identify up to isogeny the abelian subvarieties $B_W$ in terms of the group action.
In most cases we identify them in terms of Prym varieties corresponding to
subgroups $M \subset N \subset G$. We will see however that this is not the case for the dihedral group $D_{2p}$ of 
order $4p$ where $p > 3$ is a prime. We apply for this Corollary 3.10 several times 
and have to verify conditions (i) and (ii) in each case.
Note that the Prym varieties and similarly the subgroups $M \subset N \subset G$
giving the abelian subvarieties $B_i$ are not unique in general.\\

The examples,  except the quaternion group, are treated elsewhere, $S_3$ in [RR1],
$S_4$ in [RR2], $A_5$ in [SA], and $D_n$ in [CRR], but with adhoc methods, i.e. without the use
of Corollary 3.10, which was inspired by [SA].\\
\bigskip

\noindent
{\bf 4.1. The symmetric group $S_4$}

\vspace{0.5cm}
\noindent
The irreducible $\Q$-representations of $S_4$ are the trivial representation $U$, 
the alternating representation $U'$ (both of dimension 1), the standard representation 
$V$, the representation $V \otimes U'$ (both of dimension 3) and a 2-dimensional 
representation $W$ which is the pullback of the standard representation of $S_3$ 
via the canonical homomorphism $S_4 \lra S_3 = S_4/K$, where $K$ denotes the Klein subgroup.
Since every irreducible $\Q$-representation of $S_4$ is absolutely irreducible, we have 
$$
\Q [S_4] \simeq U \oplus U' \oplus W^2 \oplus V^3 \oplus (V \otimes U')^3.
$$ 
as $\Q[S_4]$-modules.
So if $X$ denotes a curve with $S_4$-action, Propositions 1.1 and 2.1 yield
$$
JX \sim A_U \times A_{U'} \times A_W \times A_{V} \times A_{V \otimes U'}
$$
\noindent
$$
A_U = B_U \,\, , \,\, A_{U'} = B_{U'} \,\, , \,\, A_W \sim B_W^2 \,\, , \,\, 
A_{V} \sim B_{V}^3 \,\, , \,\, 
A_{V \otimes U'} \sim B^3_{V \otimes U'}.
$$
We already saw in Section 3 that for $Y:= X/S_4$:
$$
A_U \sim JY.
$$
It remains to identify $A_{U'}, B_W, B_{V}$ and $ B_{V \otimes U'}$ up to isogeny 
in terms of Prym varieties.\\

\noindent
First we have
$$
\chi_{ind^{S_4}_{1_{A_4}}} = \chi_{ind^{S_4}_{1_{S_4}}} + \chi_{U'} \quad \mbox{and} \quad 
\dim_{\C} (U')^{A_4} = 1.
$$
So Corollary 3.10 implies 
$$
A_{U'} \sim P(X_{A_4}/Y).
$$
\noindent
Let $S_3$ denote the stabilizer of a symbol, say 1, in $S_4$. It is related to the standard 
representation $V$ as follows
$$
\chi_{ind^{S_4}_{1_{S_3}}} = \chi_{ind^{S_4}_{1_{S_4}}} + \chi_{V} \quad \mbox{and} \quad 
\dim_{\C} V^{S_{3}} = 1.
$$
So Corollary 3.10 gives
$$
A_V \sim P(X_{S_3}/Y)
$$
\noindent
Let $D_4$ denote the subgroup of $S_4$ generated by the cycles $(13)$ and $(1234)$, consider the subgroup 
$Z:= \,\,<(1234)>$ and compute the character table
\vspace{0.4cm}
\begin{center}
\begin{tabular}{l|c|c|c|c|c|}
& 1& (12)& (123)& (1234)& (12)(34)\\
\hline
$\chi_{ind^{S_4}_{1_Z}}$& 6&0&0&2&2\\
\hline
$\chi_{ind^{S_4}_{1_{D_4}}}$ &3&1&0&1&3\\
\hline
$\chi_{V \otimes U'}$ & 3& -1& 0&1&-1\\
\hline
\end{tabular}
\end{center}
\vspace{0.4cm} 
It is easy to check that
$$ 
\chi_{ind^{S_4}_{1_{D_4}}} = \chi_{ind^{S_4}_{1_{S_4}}} + \chi_{W}  \quad \mbox{and} \quad
\dim_{\C} W^{D_4} = 1.
$$
So Corollary 3.10 gives 
$$
B_W \sim P(X_{D_4}/Y).
$$
\noindent
Finally we deduce from the above character table
$$
\chi_{ind^{S_4}_{1_Z}} = \chi_{ind^{S_4}_{1_{D_4}}} + \chi_{V \otimes U'}.
$$
Moreover, since $\dim(V \otimes U')^Z = 1$, we may apply Corollary 3.10 to get
$$
B_{V \otimes U'} \sim P(X_{Z}/X_{D_4}).
$$
Combining everything we obtain with the notation as above\\

\noindent
{\bf Proposition 5.1:}  $JX \sim JY \times P(X_{A_4}/Y) \times P(D_4 /Y)^2 
\times P(X_{S_3}/Y)^3 \times P(X_{Z}/X_{D_4})^3$\\
\bigskip
  
\noindent
{\bf 4.2. The groups $A_5$ and $D_p$}\\
\medskip

The same method as outlined in Section 5.1 works also in the case of the alternating group $A_5$ and
the dihedral group $D_p$, for an odd prime $p$. We do not 
give the details, but only state the results.\\
 
\noindent
The alternating group $A_5$ has 5 irreducible $\C$-representations, the trivial one $U$, 
a 4-dimensional $V$ and a 5-dimensional W, all defined over $\Q$, and two 3-dimensional ones
$T_1$ and $T_2$ such that $T := T_1 \oplus T_2$ is defined over $\Q$. Hence if $X$ is a curve with 
$A_5$-action, we have
$$
JX \sim JY \oplus B_V^4 \oplus B_W^5 \oplus B_{T}^3.
$$
Proceeding as in the previous example one gets
$$
B_V \sim P(X_{A_4}/Y) \,\,\,, \,\,\, B_W \sim P(X_{D_5}/Y) \,\, \mbox{and} \,\,
 B_{T} \sim P(X_{Z_5}/X_{D_5}),
$$
where $A_4 = < (23)(45), (345)> \,\,,\,\, D_5 = <(12345), (25)(34)>$ and $Z_5 = <(12345)>$.\\

\noindent
Let $p$ be an odd prime. The dihedral group $D_p := < r,s \,\,:\,\, r^p =s^2 = (rs)^2 = 1> $ 
has two 1-dimensional representations $U$ (trivial) and $U'$ 
and $\frac{p-1}{2}$ irreducible complex representations of dimension 2:
$$
V_i (r^k) = \left( \begin{array}{cc}
\omega^{ik} & 0\\
0 & \omega^{-ik}
\end{array}
\right) , \quad V_i (sr^k) = \left( \begin{array}{cc}
0 & \omega^{-ik}\\
\omega^{ik} & 0
\end{array}
\right) \eqno(4.1)
$$
for  $i = 1, \ldots , \frac{p-1}{2}$, where $ \omega = e^{\frac{2 \pi i}{p}}$ .
There is an irreducible $\Q$-representation $W$ such that
$W \otimes_{\Q} \C := V_1 \oplus \ldots \oplus V_{\frac{p-1}{2}}.$    
So as before if $X$ is a curve with $D_p$-action we have
$$
JX \sim JY \times B_{U'} \times B_W^2.
$$
Proceeding as in the previous cases we get
$$
B_{U'} \sim P(X_{<r>}/Y) \,\,\, \mbox{and} \,\,\, B_W \sim  P(X_{<s>}/Y).
$$
\bigskip

\noindent
{\bf 4.3. The quaternion group} 
\bigskip

\noindent
For the quaternion group $Q_8 := \{ \pm 1, \pm i, \pm j, \pm ij \},
\,\,\, i^2 = j^2 = (ij)^2 = -1$ we cannot apply Corollary 3.10, 
but we can
still use  orthogonal cancellation.\\

\noindent
The center of $Q_8 $ is $< -1 >$ and
we observe that $Q_8 /< -1> \, \cong  \Z / 2 \Z \oplus \Z / 2 \Z$.
$Q_8$  has four
1-dimensional irreducible rational representations:
\vspace{0.4cm}
\begin{center}
\begin{tabular}{c|c|c|c|c|c|}
& 1& -1&$i, -i$&$j, -j$ & $ij, - ij$\\
\hline
$U$& 1&1&1&1&1\\
\hline
$U_i$& 1&1&1&-1& -1\\
\hline
$U_j$& 1&1&-1& 1&-1\\
\hline
$U_{ij}$ & 1&1&-1& -1& 1\\
\hline
\end{tabular}
\end{center}
\vspace{0.4cm}
and a $2$-dimensional complex irreducible representation $V$:
$$
i \ra \left(
\begin{array}{cc}
i & 0\\
0 & -i
\end{array}
\right) \quad , \quad j \ra 
\left(
\begin{array}{cc}
0 & i\\
i & 0
\end{array}
\right).
$$
The associated irreducible representation over $\Q$ is $W := \Ha_{\Q}$, the quaternion algebra 
\newline (-1,-1) over $\Q$
with $\Ha = \Ha_{\Q} \otimes_{\Q} \C \simeq V^2$.
So if $X$ denotes a curve with $Q_8$-action, 
Propositions 1.1 and 2.1 yield
$$
JX \sim JY \oplus B_{U_i} \oplus B_{U_j} \oplus B_{U_{ij}} \oplus B_W.
$$

\noindent
We have a diagram of coverings of curves:
$$
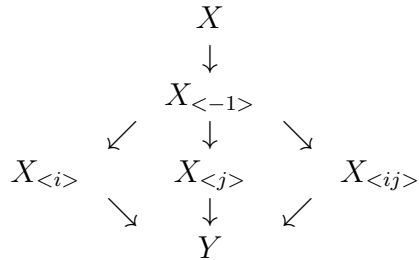
\begin{array}{ccccc}
&&X&&\\
&&\downarrow &&\\
&&X_{<-1>}&&\\
&\swarrow& \downarrow & \searrow &\\
X_{<i>} && X_{<j>} && X_{<ij>}\\
& \searrow & \downarrow & \swarrow &\\
&& Y &&
\end{array}
$$
Since $<-1>  \lhd \, Q_8$ is a normal subgroup, $Q_8 $ acts on
$X_{<-1>}$ (as $Q_8 /<-1>)$ hence on $P(X/X_{<-1>})$ and on $JX_{<-1>}$, and the 1-dimensional
irreducible rational representations of $Q_8$ are the pullbacks via $Q_8
\ra Q_8 / <-1> \,\, \cong  \,\, \Z / 2 \Z \oplus \Z / 2 \Z$ of the four
irreducible representations of $\Z / 2\Z \oplus \Z / 2 \Z$. So by 
orthogonal cancellation  the isotypical decomposition of $JX$ splits and gives:
$$
\begin{array}{rcl}
A_U \times A_{U_i} \times A_{U_j} \times A_{U_{ij}} & \sim& JX_{<-1>}\\
A_{V} & \sim  & P(X/X_{<-1>}).
\end{array}
$$
In the case of a $\Z / 2 \Z \oplus \Z / 2 \Z$-action it is easy to see (see for example [RR2]) that
$$
JX_{<-1>} \sim JY \times P(X_{<i>} / Y) \times P(X_{<j>} / Y) \times P(X_{<ij>} / Y)
$$
is the isotypical decomposition. We obtain\\

\noindent
{\bf Proposition 5.2:} {\it The isotypical decomposition of $JX$ is:}
$$
JX \sim JY \times P(X_{<i>} / Y) \times P(X_{<j>} / Y) \times P(X_{<ij>} / Y) \times P(X/X_{<-1>}).
$$

\bigskip

\noindent
{\bf 4.4. The dihedral group $D_{2p}$} 
\bigskip

Finally we give an example of a Jacobian for which not all $B'$s are Pryms of intermediate coverings.\\

Let $X$ be a curve with $D_{2p}$-action, where $D_{2p} = <r,s \,\,:\,\, r^{2p} = s^2 = (rs)^2 = 1>$
with $p > 3$ prime. Recall that $D_{2p}$ has 4   1-dimensional irreducible $\Q$-representations 
$W_1, \ldots , W_4$ and $p-1 \,\,\,$   2-dimensional irreducible complex representations $V_1, \ldots V_{p-1}$
defined as in (4.1) where now $\omega$ is a primitive $2p$-th root of unity. If we define
$$ 
W_5 = \oplus_{i=1}^{\frac{p-1}{2}} V_{2i-1}\,\,\, \mbox{and} \,\,\,
W_6 = \oplus_{i=1}^{\frac{p-1}{2}} V_{2i},
$$
 then 
$W_1, \ldots , W_6$ are the irreducible $\Q$-representations of $D_{2p}$. 
According to Propositions 1.1 and 2.1
we get
$$
JX \sim JY \times B_2 \times B_3 \times B_4 \times B_5^2 \times B_6^2
$$
with
$$
B_2 = P(X_{<r>}/Y),\,\,\, B_3 = P(X_{D_p}/Y), \,\,\, B_4 = P(X_{\tilde{D}_p} / Y) \,\,\,
\mbox{and} \,\,\, B_6 = P(X_{D_2}/Y), 
$$
where $D_p = <s,r^2>, \,\,\, \tilde{D}_p = <r^ps, r^2>$ and $D_2 = < 1,s,r^p, r^ps>$.\\

By looking at $T_0P(X_M/X_N)$ for all subgroups $M \subset N \subset D_{2p}$, one can see that none 
of those tangent spaces can be represented only in terms of $W_5$, so $B_5$ cannot be an 
intermediate Prym variety. In fact
$$
B_5 \sim (P(X_{<s>}/X_{D_p}) \cap P(X_{<s>}/X_{D_2}))^0 \subset JX_{<s>}
$$
which can also be described as the orthogonal complement of $q^*P(X_{D_s}/Y)$ in $P(X_{<s>}/X_{D_2})$ 
with respect to the induced polarization, where $q: X_{<s>} \ra X_{D_s}$. This example shows that one 
has to search other geometric objects coming from the $G$-action on the curve (see [LR]).

\begin{center}
{\bf References}
\end{center}
\begin{itemize}
\item[{[CW]}] Chevalley, C., Weil, A.: \"Uber das Verhalten der Integrale erster Gattung bei 
Automorphismen des Funktionenk\"orpers, Hamb. Abh. 10 (1934), 358 - 361.
\item[{[CRR]}] Carocca, A., Recillas, S., Rodr\'{\i}guez, R.: 
Dihedral groups acting on Jacobians, Contemp. Math. 311 (2002) 41-77.
\item[{[D]}] Donagi, R.: Decomposition of spectral covers, Ast\'{e}risque 218 (1993), 145 - 175.
\item[{[DM]}] Donagi, R., Markman, E.: Spectral covers, algebraically completely integrable Hamiltonian
systems, and moduli of bundles, in: Integrable Systems and Quantum groups, Springer LNM 1620 (1196), 1 -
119.
\item[{[H]}] Hurwitz, A.: \"Uber algebraische Gebilde mit eindeutigen
Transformationen in sich, Math. Ann. 41, (1893), 403-442.
\item[{[K]}] Kleiman, S. L.:  Relative duality for quasi-coherent sheaves, Comp. Math. 41 (1980),
39 - 60.
\item[{[LB]}] Lange, H., Birkenhake, Ch.: 
Complex Abelian Varieties.  Grundlehren No. 317, Springer (1992).
\item[{[LR]}] Lange, H., Recillas, S.: Prym Varieties of Pairs of Coverings, to appear in Adv. in Geom.
\item[{[M]}] M\'{e}rindol, J.-Y.: Vari\'{e}t\'{e}s de Prym d'un rev\^{e}tement
galoisien, 
J. reine Ang. Math. {\it 461} (1995), 49-61.
\item[{[Re]}] Recillas, S.: La Jacobiana de la extension de 
 Galois de una curva trigonal, Aport. Matem. Com. 14
(1994), 159 - 167.
\item[{[Ri]}] Ries, J.: The Prym variety for a cyclic unramified cover of a hyperelliptic curve, J.
reine ang. Math. 340 ( 1983), 59 - 69.
\item[{[RR1]}] Recillas, S., Rodr\'{\i}guez, R.: Jacobians and representations of $S_3$, Aport. Mat. Invest.
13 (1998), 117 - 140.
\item[{[RR2]}] Recillas, S., Rodr\'{\i}guez, R.: Prym varieties and four fold covers, Publ. Preliminares
Inst. Mat. Univ. Nac. Aut. Mexico, 686 (2001).
\item[{[S]}] Serre, J.-P.: R{e}pr\'{e}sentations lin\'{e}aires des groupes
finis, Herrmann (1967).
\item[{[SA]}] S\'{a}nchez-Arg\'{a}ez, A: Acciones del grupo $A_5$ en 
variedades jacobianas, Aport. Mat.
 Com.  25, (1999), 99-108.
\item[{[V]}] Viehweg, E.: Positivity of direct image sheaves and
applications to families of higher dimensional manifolds, ICTP Lect. notes 6, 249-284, Trieste (2001).
\item[{[W]}] Wirtinger, W.: Untersuchungen \"uber Theta Funktionen, Teubner, Berlin (1895).
\end{itemize}
\vspace{2cm}
\begin{center}
\begin{tabular}{ll}
Herbert Lange & Sevin Recillas\\
Mathematisches Institut & Instituto de Matematicas\\
Bismarckstr. 1 1/2 & UNAM Campus Morelia\\
D-91054 Erlangen & Morelia, Mich., 58089\\
Germany &M\'{e}xico\\
e-mail: lange@mi.uni-erlangen.de & email: sevin@matmor.unam.mx\\
&\\
&and\\
&CIMAT\\
&Callejon de Jalisco s/n\\
&Valenciana, Gto., 36000\\
&M\'{e}xico\\
&e-mail: sevin@cimat.mx
\end{tabular}
\end{center}
\end{document}